\documentclass[14pt]{amsart}
\usepackage{amsmath,amsthm,amssymb,cite}
\begin{document}
\title{VARIOUS NON-AUTONOMOUS NOTIONS FOR BOREL MEASURES}       
\author{PRAMOD DAS, TARUN DAS} 
\maketitle

\begin{abstract}
We introduce and investigate the notions of expansiveness, topological stability and persistence for Borel measures with respect to time varying bi-measurable maps on metric spaces. We prove that expansive persistent measures are topologically stable in the class of all time varying homeomorphisms.
\end{abstract}   
\medskip
   
\textbf{Mathematics Subject Classifications (2010):} Primary 54H20, Secondary 37B55     
\medskip

\textbf{Keywords:} Expansive homeomorphisms, Expansive Measures, Shadowing, Persistence, Topological Stability. 
\medskip

\textbf{(1) INTRODUCTION}
\medskip

For several decades, a discrete dynamical system induced by a continuous map or a homeomorphism on a compact metric space has been the most popular and attractive formulation for a dynamical system to a large number of mathematicians all over the world. Besides that, a significant amount of research has been carried out for group actions on compact metric spaces. Moreover, many mathematicians looked into behaviours of continuous maps or homeomorphisms on non-compact, non-metrizable spaces. 
\medskip

One of the broadly studied \cite{AH} dynamical notions in topological dynamics is expansiveness which was introduced  \cite{U} by Utz in the middle of the twentieth century. On the other hand, the most fundamental topological dynamical notions of shadowing was originated from Anosov's closing lemma \cite{A}. The \textit{Walter's Stability Theorem} \cite{W}  is one of the finest results, which states that Anosov diffeomorphisms are topologically stable. This theorem has been extended \cite{AH} to homeomorphisms on compact metric spaces using the fact that Anosov diffeomorphisms are expansive and have shadowing \cite{W1978}. In\cite{CCL}, authors have improved this result by showing that persistent (weaker notion than shadowing) expansive homeomorphisms are topologically stable. In \cite{DLRW}, the second author of the present paper with others have defined and studied the notions of shadowing, expansiveness and topological stability for homeomorphisms on uniform spaces. In particular, they have proved that Walter's stability theorem holds for homeomorphisms on certain uniform spaces. In (\cite{CL}, \cite{ALL}), authors have generalized these results to finitely generated group actions on compact metric  spaces. Recently in \cite{LM}, Lee and Morales have introduced the notions of shadowing and topological stability for Borel measures. The notion of expansiveness for Borel measures has been thoroughly investigated by Morales in \cite{M}. In \cite{LM}, they have proved a measurable version of \textit{Walter's Stability Theorem} which states that expansive measures with shadowing are topologically stable. 
\medskip

In spite of such vast literature of significant results regarding above mentioned formulations of a dynamical system, the usefulness of non-autonomous systems to understand better the topological 
\par\noindent\rule{\textwidth}{0.4pt}
\begin{center} 
Department of Mathematics, Faculty of Mathematical Sciences, 
\\
University of Delhi, Delhi-110007.
\\
Email: tarukd@gmail.com (Tarun Das), pramod.math.ju@gmail.com (Pramod Das)
\end{center} 
entropy of so called triangular map motivated Kolyada and Snoha to study such systems \cite{KS}. These type of systems arise when one studies random dynamical systems where the map used at a given time is chosen randomly from a given sequence of maps. The presence of such systems in our daily life can be felt by watching a television screen or an electronic advertisement board whose screen is divided into several units of different colour red, green, blue and so on. Because of several clear motivation discussed in the introductory paper for non-autonomous systems, authors of \cite{TD2016} have generalized the classical spectral decomposition theorem for homeomorphisms on compact metric spaces. Most importantly in \cite{TD}, authors have proved that a non-autonomous system with expansiveness and shadowing is topologically stable in the class of all time varying homeomorphisms. Our purpose is to extend this result to expansive measures with shadowing with respect to time varying homeomorphisms on relatively compact metric spaces.     
\medskip

This paper is organized as follows. In section 2, we discuss preliminaries for self-containment of the paper. In consecutive sections, we introduce and investigate expansiveness, topological stability and persistence for Borel measures with respect to time varying bi-measurable maps. We determine the size of the set of points with converging semiorbits under time varying homeomorphisms with respect to any expansive outer regular measure on separable metric spaces. Consequently, we show that every equicontinuous time varying uniform equivalence is aperiodic with respect to expansive outer regular measure. We further show that if the set of transitive points is not negligible under a measure which is persistent with respect to a time varying bi-measurable map, then each point is non-wandering. Finally, we show that on relatively compact metric spaces an expansive, persistent measure is topologically stable in the class of all time varying homeomorphisms.  
\medskip

\textbf{(2) PRELIMINARIES}
\medskip

Throughout the paper, $\mathbb{Z}$(resp. $\mathbb{N}$) denotes the set of all (resp. non-negative) integers. We consider $(X,d)$ to be any metric space unless otherwise stated. For $n\geq 1$, let $f_n:X\rightarrow X$ be a sequence of bi-measurable maps
and let $f_0:X\rightarrow X$ be the identity map. The family $F=\lbrace f_n\rbrace_{n\in\mathbb{N}}$ is called a time varying bi-measurable map on X. The inverse of $F$ is given by $F^{-1}=\lbrace f_n^{-1}\rbrace_{n\in\mathbb{N}}$. Let us denote 
\[F_n=\begin{cases} 
f_n\circ f_{n-1}\circ ...\circ f_1\circ f_0 & \textnormal {for all $n\geq 0$},
\\
f_{-n}^{-1}\circ f_{-(n-1)}^{-1}\circ ...\circ f_1^{-1}\circ f_0^{-1} & \textnormal {for all $n<0$} 
\end{cases}\]

We call $(X,F)$ an invertible non-autonomous discrete dynamical system induced by a time varying bi-measurable map $F$. 
\medskip

The dynamics of a self-homeomorphism $f$ of a metric space $X$ is a special case, where $f_n=f$ for all $n\in\mathbb{N}$. We denote
\[F_{[i,j]}=\begin{cases}
f_j\circ f_{j-1}\circ ...\circ f_{i+1}\circ f_i & \textnormal {for any $i\leq j$}
\\
$the identity map on X$ & \textnormal {for any $i>j$}
\end{cases}\] 
\[F_{[i,j]}^{-1}=\begin{cases}
f_j^{-1}\circ f_{j-1}^{-1}\circ ...\circ f_{i+1}^{-1}\circ f_i^{-1} & \textnormal {for any $i\leq j$}
\\
$the identity map on X$ & \textnormal {for any $i>j$}
\end{cases}\]

For $k\geq 1$, we define $F^k=\lbrace g_n\rbrace_{n\in\mathbb{N}}$, where $g_n=F_{[(n-1)k+1,nk]}$. The sequence $O(x_0)=\lbrace x_n\rbrace_{n\in\mathbb{Z}}$ is called the orbit of $x_0$ under $F$. Observe that $O(x_0)=\lbrace F_n(x_0)\rbrace_{n\in\mathbb{Z}}$. A subset $Y\subset X$ is said to be $F$-invariant if $f_n(Y)\subset Y$ for all $n\in\mathbb{N}$, equivalently, $F_n(Y)\subset Y$ for all $n\in\mathbb{Z}$. A homeomorphism $h:X\rightarrow X$ is called an uniform equivalence if both $h$ and $h^{-1}$ are uniformly continuous. 
\medskip

Let $(X,d_1)$ and $(Y,d_2)$ be two metric spaces. Let $F=\lbrace f_n\rbrace_{n\in\mathbb{N}}$ and $G=\lbrace g_n\rbrace_{n\in\mathbb{N}}$ be time varying bi-measurable maps on $X$ and $Y$ respectively. Then, $F$ and $G$ are said to be topologically conjugate if there is a homeomorphism $h:X\rightarrow Y$ such that $h\circ f_n=g_n\circ h$ for all $n\in\mathbb{N}$. In addition, if $h$ is an uniform equivalence, we say that $F$ and $G$ are uniformly conjugate.  
\medskip

Let $\mathcal{P}(X)$ be the power set of $X$ and let $H:X\rightarrow\mathcal{P}(X)$ be a set valued map of $X$. We define the domain of $H$ by $Dom(H)=\lbrace x\in X\mid H(x)\neq\phi\rbrace$. $H$ is said to be compact valued if $H(x)$ is compact for each $x\in X$. For some $\epsilon>0$, we write $d(H,Id)<\epsilon$ if $H(x)\subset B(x,\epsilon)$ for each $x\in X$, where $B(x,\epsilon)$ is the open $\epsilon$-ball with center $x$. Similarly, $B[x,\epsilon]$ denotes the closed $\epsilon$-ball with center $x$. $H$ is called upper semi-continuous if for every $x\in Dom(H)$ and every open neighbourhood $O$ of $H(x)$ there is $\delta>0$ such that $H(y)\subset O$ for all $y\in X$ with $d(x,y)<\delta$. A family $\mathcal{F}$ of functions is equicontinuous if for every $\epsilon>0$ there is $\delta>0$ such that $d(f(x),f(y))<\epsilon$ for all $f\in\mathcal{F}$, whenever $d(x,y)<\delta$. 
\medskip

A point $x\in X$ is called an atom for a measure $\mu$ if $\mu(x)>0$. A measure $\mu$ on $X$ is said to be non-atomic if it has no atom. Let us denote the set of all Borel measures and the set of all non-atomic Borel measures by $M(X)$ and $NAM(X)$ respectively.  
\medskip

\textbf{(3) EXPANSIVE MEASURES}
\medskip

In this section, we introduce and investigate the notion of expansiveness for Borel measures with respect to time varying bi-measurable map.
\medskip

\textbf{Definition 3.1} Let $F=\lbrace f_n\rbrace_{n\in\mathbb{N}}$ be a time varying bi-measurable map on a metric space $(X,d)$ and let $\mu\in NAM(X)$. Then, $\mu$ is said to be expansive with respect to $F$ if there is $\delta>0$ such that $\mu(\Gamma_{\delta}(x))=0$ for all $x\in X$, where $\Gamma_{\delta}(x)=\lbrace y\in X\mid d(F_n(x),F_n(y))\leq\delta$ for all $n\in\mathbb{Z}\rbrace$. Such $\delta$ is called expansive constant for $\mu$.
\medskip
 
\textbf{Remark 3.2} (i) If the space is complete separable without isolated points, then every non-atomic Borel measure is expansive with respect to any expansive \cite{TD} time varying bi-measurable map.  
\\
(ii) If $X$ is compact, then expansiveness of measure does not depend on the choice of the metric.  
\medskip 

\textbf{Example 3.3} Let $\mu$ be an expansive measure for a self-homeomorphism $f$ of a metric space $X$ and $Is$ is an isometry on $X$. Then, $\mu$ is expansive with respect to $F=\lbrace f_n\rbrace_{n\in\mathbb{N}}$, where $f_n=Is$ for $n=1,3,6,10,15,...$ and $f_n=f$ for otherwise.       
\medskip

\textbf{Proposition 3.4} Let $(X,d_1)$ and $(Y,d_2)$ be two metric spaces. Let $F=\lbrace f_n\rbrace_{n\in\mathbb{N}}$ and $G=\lbrace g_n\rbrace_{n\in\mathbb{N}}$ be time varying bi-measurable maps on $X$ and $Y$ respectively such that $F$ is uniformly conjugate to $G$. Then, $\mu$ is expansive with respect to $F$ if and only if it is expansive with respect to $G$. 
\medskip

\textbf{Proof.} Let $\mu$ be expansive with respect to $F$ with an expansive constant $\epsilon>0$. Let $h:X\rightarrow Y$ be a uniform equivalence such that $h\circ f_n=g_n\circ h$ for all $n\in\mathbb{N}$. Then, observe that $F_n\circ h^{-1}=h^{-1}\circ G_n$ for all $n\in\mathbb{N}$. Since $h^{-1}$ is uniformly continuous, there is $\delta>0$ such that $d_2(x,y)\leq\delta$ implies $d_1(h^{-1}(x),h^{-1}(y))\leq\epsilon$. Let us fix $x\in Y$.   
\begin{center}
Then, $\mu(\lbrace y\in Y\mid d_2(G_n(x),G_n(y))\leq\delta$ for all $n\in\mathbb{Z}\rbrace$
\\
$\leq \mu(\lbrace y\in Y\mid d_1(h^{-1}(G_n(x)),h^{-1}(G_n(y)))\leq\epsilon$ for all $n\in\mathbb{Z}\rbrace$
\\
$=\mu(\lbrace y\in Y\mid d_1(F_n(h^{-1}(x)),F_n(h^{-1}(y)))\leq\epsilon$ for all $n\in\mathbb{Z}\rbrace=0$.
\end{center}
 
This shows that $\mu$ is expansive with respect to $G$ with expansive constant $\delta>0$. The converse holds in similar manner because of the fact that $h$ is a uniform equivalence.
\medskip

\textbf{Theorem 3.5} Let $(X,d)$ be a metric space and let $F=\lbrace f_n\rbrace_{n\in\mathbb{N}}$ be a time varying bi-measurable map on $X$. Then, $\mu\in NAM(X)$ is expansive with respect to $F$ if and only if it is expansive with respect to $F^{-1}$.
\medskip

\textbf{Proof.} The proof is easy to work out once the fact that $F_{-n}=(F^{-1})_n$ for all $n\in\mathbb{Z}$, is clear. 
\medskip 

\textbf{Theorem 3.6} Let $F=\lbrace f_n\rbrace_{n\in\mathbb{N}}$ be a time varying uniform equivalence on a metric space $(X,d)$ such that $\lbrace f_n,f_n^{-1}\rbrace_{n\in\mathbb{N}}$ is equicontinuous. Then, $\mu\in NAM(X)$ is expansive with respect to $F$ if and only if it is expansive with respect to $G=F^k=\lbrace F_{[(n-1)k+1,nk]}\rbrace_{n\in\mathbb{N}}$ for all $k\in\mathbb{Z}$.  
\medskip

\textbf{Proof.} In view of \textit{Theorem 3.5}, it is enough to prove the result for $k\geq 1$. Let us fix $k\geq 1$ and let $e$ be an expansive constant for $\mu$. Since $\lbrace f_n,f_n^{-1}\rbrace_{n\in\mathbb{N}}$ is equicontinuous, for any $n\in\mathbb{Z}$ and any $j$ with $nk+1\leq j\leq (n+1)k$ the homeomorphisms $F_{[nk+1,j]}$ are uniformly continuous. 
\medskip

Thus, there is $\delta_j>0$ such that $d(x,y)<\delta_j$ implies $d(F_{[nk+1,j]}(x),F_{[nk+1,j]}(y))<e$ for any $n\in\mathbb{Z}$ and all $j$ with $nk+1\leq j\leq (n+1)k$. Observe that $\delta_j$ does not depend on $n$ because of the equicontinuity of $\lbrace f_n, f_n^{-1}\rbrace_{n\in\mathbb{N}}$. Then, $d(x,y)<\delta$ implies $d(F_{[nk+1,j]}(x),F_{[nk+1,j]}(y))<e$ for all $n\in\mathbb{Z}$, where $\delta=$min$\lbrace \delta_j\mid nk + 1\leq j\leq (n + 1)k\rbrace$. Observe that for any $j\in\mathbb{Z}$ there is $n\in\mathbb{Z}$ such that $nk+1\leq j\leq (n+1)k$ and $G_n=F_{nk}$ for all $n\in\mathbb{N}$ and all $k\in\mathbb{Z}$. 
\begin{center}
Now, $\mu(\lbrace y\in X\mid d(G_n(x),G_n(y))\leq \delta$ for all $n\in\mathbb{Z}\rbrace$
\\
$=\mu(\lbrace y\in X\mid d(F_{nk}(x),F_{nk}(y))\leq\delta$ for all $n\in\mathbb{Z}\rbrace$
\\
$=\mu(\lbrace y\in X\mid d(F_{[nk+1,j]}(F_{nk}(x)),F_{[nk+1,j]}(F_{nk}(y)))\leq e$ for all $n\in\mathbb{Z}$ and all $j\in\mathbb{Z}\rbrace$
\\
$=\mu(\lbrace y\in X\mid d(F_j(x),F_j(y))\leq e$ for all $j\in\mathbb{Z}\rbrace=0$.
\end{center} 

This shows that $\mu$ is expansive with respect to $F^k$ with expansive constant $\delta$.  
\medskip

Conversely, suppose that $\mu$ is expansive with respect to $F^k$ with expansive constant $\delta$. Then, for $x\in X$, $\mu(\lbrace y\in X\mid d(F^k_n(x),F^k_n(y))\leq\delta$ for all $n\in\mathbb{Z}\rbrace=0$ which implies $\mu(\lbrace y\in X\mid d(F_{nk}(x),F_{nk}(y))\leq\delta$ for all $n\in\mathbb{Z}\rbrace=0$. This further implies $\mu(\lbrace y\in X\mid d(F_n(x)F_n(y))\leq\delta$ for all $n\in\mathbb{Z}\rbrace=0$. Therefore, $\mu$ is expansive with respect to $F$ with expansive constant $\delta$.   
\medskip

\textbf{Proposition 3.7} Let $(X,d)$ be a metric space and let $F=\lbrace f_n\rbrace_{n\in\mathbb{N}}$ be a time varying bi-measurable map on $X$. Let $Y\subset X$ be an $F$-invariant subset of $X$. If $\mu\in NAM(X)$ is expansive with respect to $F$, then it is expansive with respect to the restriction of $F$ to $Y$ which is given by $F\mid_Y=\lbrace f_n\mid_Y\rbrace_{n\in\mathbb{N}}$. 
\medskip
 
\textbf{Proof.} The proof is left for the reader as an easy exercise.
\medskip

Morales have defined \cite{M} the concept of $\mu$-generators analogously as the concept of generators \cite{KR}. We define and study the concept of $\mu$-generators for invertible non-autonomous discrete dynamical system. 
\medskip

\textbf{Definition 3.8} Let $F=\lbrace f_n\rbrace_{n\in\mathbb{N}}$ be a time varying bi-measurable map on a compact metric space $(X,d)$ and let $\mu\in M(X)$ be given. Then, a finite open cover $\mathcal{A}$ of $X$ is said to be a $\mu$-generator for $F$ if for every bi-sequence $\lbrace A_n\rbrace_{n\in\mathbb{Z}}\subset\mathcal{A}$, we have $\mu(\bigcap_{n\in\mathbb{Z}} F_n(cl(A_n))=0$.
\medskip

\textbf{Theorem 3.9} Let $(X,d)$ be a compact metric space and let $F=\lbrace f_n\rbrace_{n\in\mathbb{N}}$ be a time varying bi-measurable map. Then, $\mu\in NAM(X)$ is expansive with respect to $F$ if and only if $F$ has a $\mu$-generator. 
\medskip

\textbf{Proof.} Suppose that $\mu$ is expansive with respect to $F$ with expansive constant $\delta>0$. Let $\mathcal{A}$ be the collection of all open $\delta$-balls centered at $x\in X$. Then, for any bi-sequence $\lbrace A_n\rbrace_{n\in\mathbb{Z}}\subset\mathcal{A}$, we have $\bigcap_{n\in\mathbb{Z}} F_n(cl(A_n))\subset \Gamma_{\delta}(x)$ for all $x\in \bigcap_{n\in\mathbb{Z}} F_n(cl(A_n))$. So, $\mu(\bigcap_{n\in\mathbb{Z}} F_n(cl(A_n)))\leq \mu(\Gamma_{\delta}(x))=0$. Thus, $\mathcal{A}$ is a $\mu$-generator for $F$. 
\\
Conversely, suppose that $\mathcal{A}$ is a $\mu$-generator for $F$. Let $\delta>0$ be the Lebesgue number for $\mathcal{A}$. If $x\in X$, then for every $n\in\mathbb{Z}$ there is $A_n\in \mathcal{A}$ such that $cl(B_{\delta}(F_n(x)))\subset cl(A_n)$. Thus, $\Gamma_{\delta}(x)\subset \bigcap_{n\in\mathbb{Z}} F_{-n}(cl(A_n))$, which implies $\mu(\Gamma_{\delta}(x))=0$. Therefore, $\mu$ is expansive with respect to $F$.  
\medskip

\textbf{Definition 3.10} \cite{TD} Let $(X,d)$ be a metric space and let $F=\lbrace f_n\rbrace_{n\in\mathbb{N}}$ be a time varying homeomorphism on $X$. Then, $F$ is said to be equicontinuous if $\lbrace F_{[m,n]}, F^{-1}_{[m,n]}\mid 0\leq m\leq n\rbrace$ is an equicontinuous family of functions.  
\medskip

\textbf{Theorem 3.11} Let $(X,d)$ be a Lindel{\"o}f metric space and let $F=\lbrace f_n\rbrace_{n\in\mathbb{N}}$ be an equicontinuous time varying homeomorphism. Then, there exists no expansive measure with respect to $F$.
\medskip

\textbf{Proof.} Suppose by contradiction that $\mu$ is an expansive measure with respect to $F$. Let $e$ be an expansive
constant for $\mu$. Since $F$ is equicontinuous time varying homeomorphism, the family $\lbrace F_m,F_{(-m)}\rbrace_{m\in\mathbb{N}}$ is equicontinuous. Then, for $e$ there is $\delta>0$ such that $d(x,y)<\delta$ implies $d(F_n(x),F_n(y))<e$ for all $n\in\mathbb{Z}$. Thus, $B(x,\delta)\subset \Gamma_e(x)$ and hence, $\mu(B(x,\delta))\leq \mu(\Gamma_e(x))=0$ for all $x\in X$. Now, $\lbrace B(x,\delta)\mid x\in X\rbrace$ is an open cover for $X$ and since $X$ is Lindel{\"o}f there is $\lbrace x_i\rbrace_{i\in\mathbb{N}}$ such that $\lbrace B(x,\delta)\mid i\in\mathbb{N}\rbrace$ is an open covering for $X$. So, $\mu(X)\leq\Sigma_{i\in\mathbb{N}} \mu(B(x,\delta))$, which implies $\mu(X)=0$, which is not the case.       
\medskip

\textbf{Corollary 3.12} An equicontinuous time varying homeomorphism on a complete separable metric space can not be expansive.
\medskip

\textbf{Definition 3.13} \cite{C} Let $(X,d)$ be a metric space and let $F=\lbrace f_n\rbrace_{n\in\mathbb{N}}$ be a time varying bi-measurable map on $X$. Then,
\\
(i) the $\omega$-limit set of a point $x\in X$ is given by $\omega(F,x)=\lbrace y\in X\mid$ lim$_{k\to\infty} d(F_{n_k}(x),y)=0$  for some strictly increasing sequence of integers$\rbrace$.
\\
(ii) the $\alpha$-limit set of a point $x\in X$ is given by $\alpha(F,x)=\lbrace y\in X\mid$ lim$_{k\to\infty} d(F_{n_k}(x),y)=0$ for some strictly decreasing sequence of integers$\rbrace$.  
\medskip

We say that a point $x\in X$ has converging semiorbits under $F$ if both $\alpha(F,x)$ and $\omega(G,x)$ consist of single point. The set of such points of $F$ is denoted by $A(F)$.    
\medskip

Given $x,y\in X$ and $m,n\in\mathbb{N}^+$, we define 
\begin{center}
$A(x,y,n,m)=\lbrace z\in X\mid$ max$\lbrace d(F_{-i}(z),x),d(F_i(z),y)\rbrace\leq\frac{1}{n}$ for all $i\geq m\rbrace$
\end{center}

\textbf{Lemma 3.14} Let $F=\lbrace f_n\rbrace_{n\in\mathbb{N}}$ be a time varying bi-measurable map on a separable metric space $(X,d)$. Then, there is a sequence $\lbrace x_k\rbrace_{k\in\mathbb{N}^+}\in X$ such that 
\begin{center}
$A(F)\subset \bigcap_{n\in\mathbb{N}^+}\bigcup_{k,k',m\in\mathbb{N}^{+}} A(x_k,x_{k'},n,m)$ 
\end{center}
\medskip

\textbf{Proof.} If $z\in A(F)$, then $\alpha(F,z)$ and $\omega(F,z)$ reduce to single points $x$ and $y$ respectively. Then, for each $n\in\mathbb{N}^+$ there is $m\in \mathbb{N}^+$ such that $d(F_{-i}(z),x)\leq\frac{1}{2n}$ and $d(F_i(z),y)\leq\frac{1}{2n}$ for all $i\geq m$. If $\lbrace x_k\rbrace_{k\in\mathbb{N}^+}$ is dense in $X$, there are $k,k'\in\mathbb{N}^+$ such that $d(x_k,x)\leq\frac{1}{2n}$ and $d(x_{k'},y)\leq\frac{1}{2n}$. Therefore, max$\lbrace d(F_{-i}(z),x_k),d(F_i(z),x_{k'})\rbrace\leq\frac{1}{n}$ for all $i\geq m$. This completes the proof. 
\medskip

\textbf{Lemma 3.15} Let $\mu$ be a Borel measure on a topological space. Then, for every measurable Lindel{\"o}f subset $K$ with $\mu(K)>0$ there are $z\in K$ and an open neighborhood $U$ of $z$ such that $\mu(K\cap W)>0$ for every open neighborhood $W\subset U$ of $z$.
\medskip

\textbf{Proof.} Otherwise, for every $z\in K$ there is open neighborhood $U_z\subset U$ satisfying $\mu(K\cap U_z)=0$. Since $K$ is Lindel{\"o}f, the open cover $\lbrace K\cap U_z: z\in K\rbrace$ of $K$ admits a countable sub-cover, i.e., there is a sequence $\lbrace z_l\rbrace_{l\in\mathbb{N}}$ in $K$ satisfying $K=\bigcup_{l\in\mathbb{N}}(K\cap U_{z_l})$. So, $\mu(K)\leq\sum_{l\in\mathbb{N}}\mu(K\cap U_{z_l})=0$, a contradiction.      
\medskip

\textbf{Theorem 3.16} Let $F=\lbrace f_n\rbrace_{n\in\mathbb{N}}$ be a time varying bi-measurable map on a separable metric space $(X,d)$. If $\mu$ is an expansive outer regular measure with respect to $F$, then the set $A(F)$ has measure zero with respect to $\mu$.  
\medskip

\textbf{Proof.} By contradiction, suppose there is $A\subset A(F)$ such that $\mu(A)>0$. By \textit{Lemma 3.14}, there is a sequence $\lbrace x_k\rbrace_{k\in\mathbb{N}}\in X$ such that $A(F)\subset \bigcap_{n\in\mathbb{N}^+}\bigcup_{k,k',m\in\mathbb{N}^{+}} A(x_k,x_{k'},n,m)$. It follows that $A\subset \bigcup_{k,k',m\in\mathbb{N}^{+}} A(x_k,x_{k'},n,m)$ for all $n\in\mathbb{N}^+$. Thus, we can choose $k,k',n,m\in\mathbb{N}^+$ with $\frac{1}{n}\leq\frac{e}{2}$ such that $\mu(A(x_k,x_{k'},n,m)>0$. Hereafter, we fix such $k,k',n,m\in\mathbb{N}^+$ and for simplicity we put $B=A(x_k,x_{k'},n,m)$. 
\medskip

Since $X$ is separable, it is a second countable metric space. Since $\mu$ is outer regular, the \textit{Lusin theorem} implies that for every $\epsilon>0$ there is a measurable set $C_{\epsilon}\subset X$ with $\mu(X\setminus C_{\epsilon})<\epsilon$ such that $F_i\mid_{C_{\epsilon}}$ is continuous for all $\mid i\mid\leq m$. Taking $\epsilon=\frac{\mu(B)}{2}$, we get a measurable set $C=C_{\frac{\mu(B)}{2}}$ such that $F_i\mid_C$ is continuous for all $\mid i\mid\leq m$ and $\mu(B\cap C)>0$. Further, since $K=B\cap C$ is a Lindel{\"o}f subspace of $X$, by \textit{Lemma 3.15} there are $z\in B\cap C$ and $\delta_0>0$ such that $\mu(B\cap C\cap B[z,\delta])>0$ for all $0<\delta<\delta_0$. Since $z\in C$ and $F_i\mid_C$ is continuous for all $\mid i\mid\leq m$, we can fix $0<\delta<\delta_0$ such that $d(F_i(z),F_i(w))\leq e$ for all $\mid i\mid\leq m$, whenever $d(z,w)\leq\delta$ with $w\in C$. 
\medskip

Now let $w\in B\cap C\cap B[z,\delta]$ which implies $w\in C\cap B[z,\delta]$ and hence, $d(F_i(z),F_i(w))\leq e$ for all $\mid i\mid\leq m$. Again, $z,w\in B=A(x_k,x_{k'},n,m)$, so observe that $d(F_i(z),F_i(w))\leq e$ for all $\mid i\mid\geq m$. Combining we get $d(F_i(z),F_i(w))\leq e$ for all $i\in\mathbb{Z}$ which implies $w\in\Gamma_e(z)$ and hence, $B\cap C\cap B[z,\delta]\subset\Gamma_e(z)$. Thus $\mu(B\cap C\cap B[z,\delta])=0$, a contradiction.     
\medskip

\textbf{Example 3.17} For $n\geq 0$, let $f_n:\mathbb{R}\rightarrow \mathbb{R}$ be given by $f_n(x)=x$ if $x\in \mathbb{Q}$ and $f_n(x)=(n+1)x$ if $x\in \mathbb{R}\setminus \mathbb{Q}$. Then, the Lebesgue measure on $\mathbb{R}$ is expansive with respect to the time varying bi-measurable map $F=\lbrace f_n\rbrace_{n\in\mathbb{N}}$. By \textit{Theorem 3.16}, $A(F)$ has measure zero with respect to $\mu$. 
\medskip

\textbf{Definition 3.18} Let $F=\lbrace f_n\rbrace_{n\in\mathbb{N}}$ be a time varying bi-measurable map on a metric space $(X,d)$. A point $p\in X$ is said to be periodic if there is an integer $k>0$ such that $F_{ik+j}(p)=F_j(p)$ for all $i\in\mathbb{Z}$ and $0\leq j<k$. The positive integer $k$ is said to be a period of $p$. A set $A\subset X$ is said to be periodic if there is $k>0$ such that each point in $A$ is periodic with period $k$. A time varying bi-measurable map $F$ is said to be aperiodic with respect to a measure $\mu$ if every measurable periodic subset of $X$ has measure zero with respect to $\mu$.     
\medskip

\textbf{Corollary 3.19} If $F=\lbrace f_n\rbrace_{n\in\mathbb{N}}$ is a time varying uniform equivalence on a separable metric space such that $\lbrace f_n,f_n^{-1}\rbrace_{n\in\mathbb{N}}$ is equicontinuous, then it is aperiodic with respect to an expansive outer regular measure.  
\medskip

\textbf{Proof.} Let $\mu$ be an expansive outer regular measure with respect to $F$. Let $m$ be a positive integer and let $A$ be a measurable subset such that for each $x\in A$, we have $F_{im+j}(x)=F_j(x)$ for all $i\in\mathbb{Z}$ and $0\leq j<m$. Thus, the orbit of $x$ is given by $\lbrace F_0(x),F_1(x),...,F_{m-1}(x)\rbrace$. Since $F_{nm}(x)=\lbrace x\rbrace$ for all $n\geq 1$, we have $A\subset A(F^m)$, where $F^m=\lbrace g_n=F_{[(n-1)m+1,nm]}\rbrace_{n\in\mathbb{N}}$. By \textit{Theorem 3.6} $\mu$ is expansive measure for $F^m$. So, by \textit{Theorem 3.16} $\mu(A)\leq \mu(A(F^m))=0$.                             
\medskip

\textbf{(4) TOPOLOGICALLY STABLE MEASURES} 
\medskip

Let $(X,d)$ be a metric space and $d_1(x,y)=$min$\lbrace d(x,y),1\rbrace$ is the standard bounded metric on $X$. Let $\mathcal{H}(X)$ be the metric space of all bi-measurable maps with the metric $\eta(f,g)=$sup$_{x\in X} d_1(f(x),g(x))$. If $\mathcal{G}(X)$ is the collection of all time varying bi-measurable maps, then we define a metric $p$ on $\mathcal{G}(X)$ as $p(F,G)=$max$\lbrace$ sup$_{n\in\mathbb{N}}$ $\eta(f_n,g_n)$, sup$_{n\in\mathbb{N}}$ $\eta(f_n^{-1},g_n^{-1})\rbrace$, where $F=\lbrace f_n\rbrace_{n\in\mathbb{N}}$ and $G=\lbrace g_n\rbrace_{n\in\mathbb{N}}$.  
\medskip

\textbf{Definition 4.1} Let $F=\lbrace f_n\rbrace_{n\in\mathbb{N}}$ be a time varying bi-measurable map on a metric space $(X,d)$. 
\\
(a) A Borel measure $\mu$ is said to be topologically stable with respect to $F$ if for every $\epsilon>0$ there is $0<\delta<1$ such that if $G=\lbrace g_n\rbrace_{n\in\mathbb{N}}$ is another time varying bi-measurable map on $X$ with $p(F,G)<\delta$, then there is an upper semi-continuous compact valued map $H:X\rightarrow \mathcal{P}(X)$ with measurable domain satisfying (i) $\mu(X\setminus Dom(H))=0$, (ii) $\mu\circ H=0$, (iii) $d(H,Id)<\epsilon$, (iv) $F_n(H(x))\subset B[G_n(x),\epsilon]$ for all $n\in\mathbb{Z}$.
\\
(b) $F$ is said to be topologically stable \cite{TD} if for every $\epsilon>0$ there is $0<\delta<1$ such that if $G=\lbrace g_n\rbrace_{n\in\mathbb{N}}$ is another time varying bi-measurable map on $X$ with $p(F,G)<\delta$, then there is a continuous map $h:X\rightarrow X$ such that $d(h(x),x)<\epsilon$ and $d(F_n(h(x)),G_n(x))<\epsilon$ for all $n\in\mathbb{Z}$.  
\medskip

\textbf{Theorem 4.2} Let $F=\lbrace f_n\rbrace_{n\in\mathbb{N}}$ be a time varying bi-measurable map with topological stability on a metric space $(X,d)$. Then, every non-atomic Borel measure $\mu$ on $X$ is topologically stable with respect to $F$.   
\medskip

\textbf{Proof.} Let us fix $\epsilon>0$ and let $0<\delta<1$ be given for $\epsilon$ by the topological stability of $F$. Let $G=\lbrace g_n\rbrace_{n\in\mathbb{N}}$ be another time varying bi-measurable map such that $p(F,G)<\delta$. Then, there is a continuous map $h:X\rightarrow X$ such that $d(h(x),x)<\epsilon$ and $d(F_n(h(x)),G_n(x))<\epsilon$ for all $n\in\mathbb{Z}$. Define the compact valued map $H:X\rightarrow\mathcal{P}(X)$ by $H(x)=\lbrace h(x)\rbrace$ for all $x\in X$. $H$ is upper semi-continuous because of continuity of $h$. Since $Dom(H)=X$, $\mu(X\setminus Dom(H))=0$ and since $\mu$ is non-atomic $\mu\circ H=0$. Further, $d(H(x),x)=d(h(x),x)<\epsilon$ which gives $d(H,Id)<\epsilon$. Finally, $F_n(H(x))=F_n(h(x))\subset B[G_n(x),\epsilon]$ for all $n\in\mathbb{Z}$.  
\medskip

\textbf{Corollary 4.3} Every complete separable metric space supporting topologically stable bi-measurable map without topologically stable measure is countable.  
\medskip

\textbf{Proof.} It follows from the well known fact that every uncountable complete separable metric space admits non-atomic Borel probability measure.     
\medskip

Given a continuous map $h:X\rightarrow X$ and a Borel measure $\mu$ on $X$, we define the measure $h_*(\mu)(A)=\mu(h(A))$ for measurable $A\subset X$.    
\medskip

\textbf{Theorem 4.4} Let $F=\lbrace f_n\rbrace_{n\in\mathbb{N}}$ be a time varying bi-measurable map on a metric space $(X,d)$ and let $h:X\rightarrow X$ be an uniform equivalence. If $\mu\in M(X)$ is topologically stable with respect to $F$, then $h_*(\mu)$ is topologically stable with respect to $F'=\lbrace f_n'\rbrace_{n\in\mathbb{N}}$, where $f_n'=h^{-1}\circ f_n\circ h$ for all $n\in\mathbb{N}$.
\medskip

\textbf{Proof.} Fix $0<\epsilon<1$. There is $0<\epsilon'<\epsilon$ such that $d_1(a,b)<\epsilon'$ implies $d_1(h^{-1}(a),h^{-1}(b))<\epsilon$. Let $\delta'>0$ be given for $\epsilon'$ by the topological stability of $\mu$ with respect to $F$. Again, there is $\delta>0$ such that $d_1(a,b)<\delta$ implies $d_1(h(a),h(b))<\delta'$. 
\medskip

Suppose $G=\lbrace g_n\rbrace_{n\in\mathbb{N}}$ be another time varying bi-measurable map such that $p(F',G)<\delta$ i.e.; max$\lbrace$sup$_{n\geq 0}$ $\eta(f_n',g_n)$, sup$_{n\geq 0}$ $\eta(f_n'^{-1},g_n^{-1})\rbrace<\delta$ which implies $d_1((h^{-1}\circ\ f_n\circ h)(x),g_n(x))<\delta$ and $d_1((h^{-1}\circ f_n^{-1}\circ h)(x),g_n^{-1}(x))<\delta$ for all $x\in X$, $n\in\mathbb{N}$. Thus, $d_1(f_n(h(x)),(h\circ g_n\circ h^{-1})(h(x)))<\delta'$ and $d_1(f_n^{-1}(h(x)),(h\circ g_n\circ h^{-1})^{-1}(h(x))<\delta'$ for all $x\in X$, $n\in\mathbb{N}$, which implies $d_1((f_n(h(x)),g_n'(h(x))<\delta'$, $d_1(f_n^{-1}(h(x)),g_n'^{-1}(h(x))<\delta'$ for all $x\in X$ and $n\in\mathbb{N}$, where $g_n'=h\circ g_n\circ h^{-1}$ for all $n\in\mathbb{N}$. Let us put $G'=\lbrace g_n'\rbrace_{n\in\mathbb{N}}$. Therefore, by topological stability of $\mu$ with respect to $F$, there is an upper semi-continuous compact valued map $H:X\rightarrow\mathcal{P}(X)$ with measurable domain satisfying $\mu(X\setminus Dom(H))=0$, $\mu\circ H=0$, $d(H,Id)<\epsilon'$, $F_n(H(x))\subset B[G_n'(x),\epsilon']$ for all $n\in\mathbb{Z}$. 
\medskip

One can verify that, $K=h^{-1}\circ H\circ h$ is an upper semi-continuous compact valued map of $X$ with measurable domain such that $(h_*\mu)(X\setminus Dom(K))=0$, $(h_*\mu)\circ K=0$, $d(K,Id)<\epsilon$. Finally, observe that for each $n\in\mathbb{Z}$, we have $F_n'(K(x))=(h^{-1}\circ F_n\circ h)(K(x))=(h^{-1}\circ F_n)(H(h(x))\subset h^{-1}(B[(h\circ G_n\circ h^{-1})(h(x)),\epsilon'])=h^{-1}(B[h(G_n(x)),\epsilon'])\subset B[G_n(x),\epsilon]$.  
\medskip

\textbf{Remark 4.5} In view of the above theorem we can conclude that if $F$ and $G$ are topologically conjugated bi-measurable maps on a metric space, then there is a bijective correspondence between the set of all topologically stable measures with respect to $F$ and the set of all topologically stable measures with respect to $G$.   
\medskip

\textbf{Theorem 4.6} Every topologically stable measure of an expansive time varying bi-measurable map is non-atomic (hence, expansive).  
\medskip

\textbf{Proof.} Let $\mu$ be a topologically stable measure with respect to an expansive time varying bi-measurable map $F=\lbrace f_n\rbrace_{n\in\mathbb{N}}$ on a metric space $(X,d)$. Let $\epsilon>0$ be an expansive constant for $F$ and let $0<\delta<1$ be given for $\epsilon$ by the topological stability of $\mu$. Taking $G=F$ in the definition of topological stability of $\mu$, we get an upper semi-continuous compact valued map $H:X\rightarrow{P}(X)$ with measurable domain such that $\mu(X\setminus Dom(H))=0$, $\mu\circ H=0$, $d(H,Id)<\epsilon$ and $F_n(H(x))\subset B[F_n(x),\epsilon]$ for all $n\in\mathbb{Z}$. If $x\in Dom(H)$, then there is $y\in H(x)$. Thus, $F_n(y)\in B[F_n(x),\epsilon]$ for all $n\in\mathbb{Z}$ and hence $d(F_n(x),F_n(y))\leq\epsilon$ for all $n\in\mathbb{Z}$. Since, $\epsilon$ is an expansive constant we must have $x=y$. Therefore, $H(x)=\lbrace x\rbrace$ for all $x\in Dom(H)$. If possible, suppose that $z$ is an atom for $\mu$. Since $\mu(X\setminus Dom(H))=0$, $z\in Dom(H)$. Thus, $H(z)=\lbrace z\rbrace$ and hence, $\mu(H(z))=0$. This is a contradiction.    
\medskip

\textbf{(5) PERSISTENT MEASURES}   
\medskip

In this section, our purpose is to study the notion of persistence for Borel measures with respect to time varying bi-measurable maps. 
\medskip
 
\textbf{Definition 5.1} Let $F=\lbrace f_n\rbrace_{n\in\mathbb{N}}$ be a time varying bi-measurable map on a metric space $(X,d)$. 
\\
(i) A Borel measure $\mu$ is said to be persistent with respect to $F$ if for every $\epsilon>0$ there is $0<\delta<1$ and a measurable set $B\subset X$ with $\mu(X\setminus B)=0$ such that if $G=\lbrace g_n\rbrace_{n\in\mathbb{N}}$ is another time varying bi-measurable map with $p(F,G)<\delta$, then for each $x\in B$ there is $y\in X$ such that $d(F_n(y),G_n(x))<\epsilon$ for all $n\in\mathbb{Z}$. 
\\
(ii) $F$ is said to be persistent if for every $\epsilon>0$ there is $0<\delta<1$ such that if $G=\lbrace g_n\rbrace_{n\in\mathbb{N}}$ is another time varying bi-measurable map with $p(F,G)<\delta$, then for each $x\in X$ there is $y\in X$ such that $d(F_n(y),G_n(x))<\epsilon$ for all $n\in\mathbb{Z}$.
\\
(iii) \cite{TD} A sequence $\lbrace x_n\rbrace_{n\in\mathbb{Z}}$ is said to be a $\delta$-pseudo orbit if $d(f_{n+1}(x_n),x_{n+1})<\delta$ for all $n\geq 0$ and $d(f_{-n}^{-1}(x_{n+1}),x_n)<\delta$ for all $n\leq -1$. A $\delta$-pseudo orbit is said to be through $B\subset X$ if $x_0\in B$. A sequence $\lbrace x_n\rbrace_{n\in\mathbb{N}}$ is said to be $\epsilon$-shadowed by some point $y$ in $X$ if $d(F_n(y),x_n)<\epsilon$ for all $n\in\mathbb{Z}$. A measure $\mu$ is said to have shadowing with respect to $F$ if for every $\epsilon>0$ there is $\delta>0$ and a measurable $B\subset X$ with $\mu(X\setminus B)=0$ such that every $\delta$-pseudo orbit through $B$ is $\epsilon$-shadowed by some point in $X$. If $B=X$, then we say that $F$ has shadowing.
\medskip

\textbf{Lemma 5.2} Let $F=\lbrace f_n\rbrace_{n\in\mathbb{N}}$ be a time varying bi-measurable map on a metric space $(X,d)$. If $\mu\in M(X)$ has shadowing with respect to $F$, then it is persistent with respect to $F$. 
\medskip

\textbf{Proof.} Let $\epsilon>0$ be given and let $\delta>0$ and measurable $B\subset X$ with $\mu(X\setminus B)=0$ be given for $\epsilon$ by the shadowing of $\mu$ with respect to $F$. Let $G=\lbrace g_n\rbrace_{n\in\mathbb{N}}$ be another time varying bi-measurable map such that $p(F,G)<\delta$. Then, one can prove that for any $x\in B$, $\lbrace G_n(x)\rbrace_{n\in\mathbb{N}}$ is a $\delta$-pseudo orbit through $B$. So by shadowing of $\mu$, there is $y\in X$ such that $d(F_n(y),G_n(x))<\delta$ for all $n\in\mathbb{Z}$.      
\medskip
 
\textbf{Remark 5.3} (i) If $F$ is persistent, then any non-atomic measure is persistent with respect to $F$.      
\\
(ii) If $X$ is compact, then persistence of measure does not depend on the choice of the metric. 
\medskip

\textbf{Theorem 5.4} \cite{TD} Let $(X,d)$ be a compact metric space and let $F=\lbrace f_n\rbrace_{n\in\mathbb{N}}$ be a time varying homeomorphism on $X$ such that $\lbrace f_n,f_n^{-1}\rbrace_{n\in\mathbb{N}}$ is an equicontinuous family. Then, $F$ has shadowing if and only if $F^k$ has shadowing for all $k\in\mathbb{Z}\setminus\lbrace 0\rbrace$.             
\medskip

\textbf{Example 5.5} Let $F=\lbrace f_n\rbrace_{n\in\mathbb{N}}$ be a time varying homeomorphism on a compact metric space $(X,d)$ which is not totally disconnected and let $f:X\rightarrow X$ be $N$-expansive homeomorphism with shadowing \cite{CC}. If $f_0=I,f_n=f$ for $n$ odd, $f_n=f^{-1}$ for $n$ even, then $F^2=\lbrace g_n\rbrace_{n\in\mathbb{N}}$ where $g_n=I$ for all $n\in\mathbb{N}$. It is well known that $F^2$ does not have shadowing. Therefore, by \textit{Theorem 5.4} $F$ can not have shadowing. This shows that $F$ may not have shadowing in spite of the fact that each $f_n$ has shadowing.  
\medskip

\textbf{Theorem 5.6} Let $F=\lbrace f_n\rbrace_{n\in\mathbb{N}}$ be a time varying bi-measurable map on a metric space $(X,d)$ and let $h:X\rightarrow X$ be a uniform equivalence. If $\mu\in M(X)$ is persistent with respect to $F$, then $h_*(\mu)$ is persistent with respect to $F'=\lbrace f_n'\rbrace_{n\in\mathbb{N}}$, where $f_n'=h^{-1}\circ f_n\circ h$ for all $n\in\mathbb{N}$. 
\medskip

\textbf{Proof.} The proof is similar to that of \textit{Theorem 4.4} above.
\medskip

\textbf{Remark 5.7} In view of the above theorem one can conclude that if $F$ and $G$ are topologically conjugated bi-measurable maps, then there is a bijective correspondence between the set of all persistent measures with respect to $F$ and the set of all persistent measures with respect to $G$.          
\medskip

\textbf{Definition 5.8} Let $F=\lbrace f_n\rbrace_{n\in\mathbb{N}}$ be a time varying bi-measurable map on a metric space $(X,d)$.
\\
(i) A point $x\in X$ is said to be non-wandering if for any non-empty open set $U$ containing $x$ and for any $n\geq 0$ there is $m\geq n$ and $r\geq 0$ such that $F_{[m,m+r]}(U)\cap U\neq\phi$ or $F_{[m,m+r]}^{-1}(U)\cap U\neq\phi$. The set of such points is denoted by $\Omega(F)$.
\\
(ii) $F$ is said to be transitive with respect to $\mu\in M(X)$ if the set of transitive points (i.e; points $x$ for which $\omega(F,x)=X$) has positive measure with respect to $\mu$.    
\medskip

\textbf{Lemma 5.9} \cite{TD2013} Let $F=\lbrace f_n\rbrace_{n\in\mathbb{N}}$ be a time varying bi-measurable map on a metric space $(X,d)$. Then, $\Omega(F)$ is a closed subset of $X$. 
\medskip 

\textbf{Proof.} Let $\lbrace x_n\rbrace_{n\in\mathbb{N}}$ be a sequence in $\Omega(F)$ converging to $x\in X$. Let $\delta>0$ and $s\in\mathbb{N}$. Since $x_n\rightarrow x$, there is $l\in\mathbb{N}$ such that $x_l\in B(\delta,x)$. Then, there is $\epsilon>0$ such that $B(\epsilon,x_l)\subset B(\delta,x)$. Now, $x_l\in \Omega(F)$ implies that there is $m\geq s$ and $r\geq 0$ such that  $F_{[m,m+r]}(B(\epsilon,x_l))\cap B(\epsilon,x_l)\neq\phi$ or $F_{[m,m+r]}^{-1}(B(\epsilon,x_l))\cap B(\epsilon,x_l)\neq\phi$. This further implies $F_{[m,m+r]}(B(\delta,x))\cap B(\delta,x)\neq\phi$ or $F_{[m,m+r]}^{-1}(B(\delta,x))\cap B(\delta,x)\neq\phi$ and hence, $x\in \Omega(F)$ shows that $\Omega(F)$ is closed.   
\medskip

\textbf{Theorem 5.10} Let $\mu\in M(X)$ be a persistent measure with respect to a time varying bi-measurable map $F=\lbrace f_n\rbrace_{n\in\mathbb{N}}$ on a relatively compact metric space $(X,d)$. If $F$ can be approximated by another time varying bi-measurable map $G=\lbrace g_n\rbrace_{n\in\mathbb{N}}$ which is transitive with respect to $\mu$, then $\Omega(F)=X$.   
\medskip

\textbf{Proof.} Let us fix $\epsilon>0$ and $z\in X$. By relative compactness of $X$, the closed bounded set $B[K,\epsilon]$ is compact for any compact subset $K$ of $X$. Let $\delta>0$ and a measurable set $B\subset X$ with $\mu(X\setminus B)=0$ be given for $\frac{\epsilon}{2}$ by the persistence of $\mu$ with respect to $F$. From the hypothesis, we have $p(F,G)<\delta$ and since $G$ is transitive with respect to $\mu$, there is $x\in B$ such that $\omega(G,x)=X$. Since $x\in B$, by persistence of $\mu$ there is $y\in X$ such that $d(F_n(y),G_n(x))<\frac{\epsilon}{2}$ for all $n\in\mathbb{Z}$. Since $\omega(G,x)=X$, for any $z\in X$, there is a sequence $n_k$ such that $d(G_{n_k}(x),z)\to 0$ as $k\to \infty$. Thus, there is $k'\in\mathbb{N}$ such that $d(G_{n_k}(x),z)<\frac{\epsilon}{2}$ for all $k\geq k'$. Therefore, $d(F_{n_k}(y),z)\leq\epsilon$ for all $k\geq k'$ which implies $F_{n_k}(y)\in B[z,\epsilon]$ for all $k\geq k'$. Since $B[z,\epsilon]$ is compact, $F_{n_k}(y)\to w$ for some $w\in X$. One can verify that $w\in \Omega(F)$. Since $d(z,F_{n_k}(y))\leq d(z,G_{n_k}(x))+d(G_{n_k}(x),F_{n_k}(y))\leq d(z,G_{n_k}(x))+\frac{\epsilon}{2}$ for all $k\in\mathbb{N}$. By letting $k\rightarrow\infty$, we get $d(z,w)<\epsilon$ which implies $d(z,\Omega(F))<\epsilon$. Thus, $\Omega(F)$ is dense in $X$. By \textit{Lemma 5.9} $\Omega(F)$ is closed, so we must have $\Omega(F)=X$.            
\medskip

\textbf{Theorem 5.11} Every topologically stable measure with respect to a time varying bi-measurable map is persistent.
\medskip

\textbf{Proof.} Let $\epsilon>0$ be fixed and $\delta>0$ be given for $\epsilon$ by the topological stability of $\mu$ with respect to a time varying bi-measurable map $F=\lbrace f_n\rbrace_{n\in\mathbb{N}}$. Let $G=\lbrace g_n\rbrace_{n\in\mathbb{N}}$ be another time varying bi-measurable map such that $p(F,G)<\delta$. Then, by topological stability of $\mu$ there is an upper semi-continuous compact valued map $H:X\rightarrow\mathcal{P}(X)$ with measurable domain such that $\mu(X\setminus Dom(H))=0=\mu\circ H$, $d(H,Id)<\epsilon$ and $F_n(H(x))\subset B[G_n(x),\epsilon]$ for all $n\in\mathbb{Z}$. If $B=Dom(H)$, then for each $x\in B$, there is $y\in H(x)$. Therefore, $F_n(y)\in B[G_n(x),\epsilon]$ for all $n\in\mathbb{Z}$ which implies $d(F_n(y),G_n(x))\leq\epsilon$ for all $n\in\mathbb{Z}$. This completes the proof. 
\medskip

In view of the above theorem, we find sufficient conditions for persistent measure to be topologically stable in the following theorem which is extension of Walter's stability theorem for homeomorphisms.  
\medskip

\textbf{Theorem 5.12} Let $(X,d)$ be a relatively compact metric space and let $F=\lbrace f_n\rbrace_{n\in\mathbb{N}}$ be a time varying homeomorphism on $X$. If $\mu\in NAM(X)$ is expansive and persistent with respect to $F$, then it is topologically stable with respect to $F$. 
\medskip

\textbf{Proof.} Let $\mu$ be expansive and persistent with respect to $F$. Let $e$ be an expansive constant for $\mu$. Take $0<\epsilon<1$ and $0<\epsilon'<$min$(\frac{e}{2},\epsilon)$. By relative compactness of $X$, the closed bounded set $B[K,\epsilon']$ is compact for any compact subset $K\subset X$. Let $0<\delta<1$ and a measurable set $B\subset X$ with $\mu(X\setminus B)=0$ be given for $\epsilon'$ by the persistence of $\mu$. Let $G=\lbrace g_n\rbrace_{n\in\mathbb{N}}$ be another time varying homeomorphism such that $p(F,G)<\delta$. Define the compact valued map $H:X\rightarrow \mathcal{P}(X)$ given by $H(x)=\bigcap_{n\in\mathbb{Z}} F_n^{-1}(B[G_n(x),\epsilon'])$  
\medskip

First we prove that $Dom(H)$ is measurable. Take a sequence $x_k\in Dom(H)$ converging to some $x\in X$. Since $x_k\in Dom(H)$, we can choose a sequence $y_k\in X$ such that $d(F_n(y_k),G_n(x_k))\leq\epsilon'$ for all $k\in\mathbb{N}$, $n\in\mathbb{Z}$. Thus, $y_k\in B[x_k,\epsilon']$ for all $k\in\mathbb{N}$. Let $K$ be a compact neighbourhood of $x$. Since $x_k\to x$ as $k\to\infty$, there is $N\in\mathbb{N}$ such that $x_k\in K$ for all $k\geq N$. So, $y_k\in B[K,\epsilon']$ for all $k\geq N$. Since $B[K,\epsilon']$ is compact, we can assume that $y_k$ converges to some point $y$ in $X$. Therefore, $d(F_n(y),G_n(x))\leq\epsilon'$ for all $n\in\mathbb{Z}$ which implies $y\in H(x)$. This shows that $x\in Dom(H)$ which means $Dom(H)$ is closed and hence, measurable. 
\medskip

We now prove that $\mu(X\setminus Dom(H))=0$. By persistence of $\mu$, for each $x\in B$ there is $y\in X$ such that $d(F_n(y),G_n(x))\leq\epsilon'$ for all $n\in\mathbb{Z}$ which implies $y\in H(x)$. This means $H(x)\neq\phi$ for all $x\in B$ and thus, $B\subset Dom(H)$. Therefore, $\mu(X\setminus Dom(H))=0$. 
\medskip

Afterwords, we prove that $H$ is upper semi-continuous. Let $x\in Dom(H)$ and $O$ be an open neighbourhood of $H(x)$. Define
$H(y)=\bigcap_{m\geq 0} H_m(y)$, where $H_m(y)=\bigcap_{n=-m}^m F_n^{-1}(B[G_n(y),\epsilon'])$. Clearly, there is $m\in\mathbb{Z}$ such that $H_m(y)\subset O$. We assert that there is $\eta>0$ such that $H_m(y)\subset O$ for all $y\in X$ with $d(x,y)<\eta$. If not, there exists $y_k\to x$ as $k\to\infty$ and $z_k\in H_m(y_k)\setminus O$ for all $k\in\mathbb{N}$. As earlier, we can assume that $z_k$ converges to a point, say $z$ and observe that $z\notin O$. But $z_k\in H_m(y_k)$ for all $k\in\mathbb{N}$ implies that $d(F_n(z_k),G_n(y_k))\leq\epsilon'$ for all $k\in\mathbb{N}$ and $-m\leq n\leq m$. Then, $d(F_n(z),G_n(x))\leq\epsilon'$ for $-m\leq n\leq m$. So, $z\in H_m(x)\subset O$. This leads to a contradiction. Hence, $H$ is upper semi-continuous. 
\medskip

Now, we prove $\mu\circ H=0$. Take $x\in X$ and $y\in H(x)$. If $z\in H(x)$, then we have $d(F_n(z),G_n(x))\leq\epsilon'$ for every $n\in\mathbb{Z}$. Since $y\in H(x)$, we have $d(F_n(y),G_n(x))\leq\epsilon'$ for every $n\in\mathbb{Z}$. Therefore,
$d(F_n(y),F_n(z))\leq 2\epsilon'$ for all $n\in\mathbb{Z}$. Since $2\epsilon'<e$, we conclude that $z\in \Gamma_e(y)$ which implies that $H(x)\subset\Gamma_e(y)$. Since $e$ is an expansive constant of $\mu$, $\mu(H(x))\leq\mu(\Gamma_e(y))=0$. 
\medskip

It follows from the definition of $H$ that $H(x)\subset B[x,\epsilon']$. Since $\epsilon'<\epsilon$, we also have $d(H,Id)\leq\epsilon$. 
\medskip

Finally, from the definition of $H(x)$ it is clear that $F_n(H(x))\subset B[G_n(x),\epsilon]$ for all $n\in\mathbb{Z}$.
\medskip 

The following is an immediate consequence of the fact that Borel measures with shadowing are persistent.  
\medskip

\textbf{Corollary 5.13} Let $(X,d)$ be a relatively compact metric space and let $F=\lbrace f_n\rbrace_{n\in\mathbb{N}}$ be a time varying homeomorphism on $X$. If $\mu\in NAM(X)$ is expansive and has shadowing with respect to $F$, then it is topologically stable with respect to $F$.  
\medskip

\textbf{Conflict of interests:} No conflict of interests is reported by the authors.     
\medskip

\textbf{Acknowledgement:} The first author is supported by Department of Science and Technology, Government of India, under INSPIRE Fellowship (Resgistration No- IF150210) since march 2015.          

\end{document}